\documentclass[letterpaper,11pt]{article}
\usepackage[utf8]{inputenc}
\usepackage{authblk}
\usepackage{amssymb}
\usepackage{amsthm}
\usepackage{amsmath}
\usepackage{amsxtra}
\usepackage{amsfonts}
\usepackage{mathtools}
\usepackage{bm}
\usepackage{siunitx}
\usepackage{tikz}
\usepackage[top=1in,bottom=1in,right=1in,left=1in]{geometry}
\usepackage{geometry}
\usepackage{tkz-euclide}
\usepackage{listings}
\usepackage{amsthm}
\usepackage{algorithm}
\usepackage{pstricks-add}
\usepackage{colonequals}
\usepackage{graphicx} 
\usepackage{epstopdf}  
\usepackage{booktabs}
\usepackage{fancyhdr}
\usepackage{algorithm}
\usepackage{algpseudocode}
\usepackage{enumitem}
\usepackage{diagbox}
\usepackage[all]{xy}
\usepackage[normalem]{ulem}
\usepackage{amsmath}
\usepackage{makecell}

\usepackage{color}
\usepackage{subfig}
\usepackage{placeins}
\usepackage{slashbox}
\usepackage{float}
\usepackage{colortbl}
\usepackage{soul}

\newtheorem{theorem*}{Theorem}

\definecolor{Gray}{gray}{0.9}

\newcommand{\yl}[1]{{\color{black}#1}}
\newcommand{\rv}[1]{{\color{black}#1}}

\title{Sampling error mitigation through spectrum smoothing:\\first experiments with ensemble transform Kalman filters and Lorenz models}
\author[1]{Bosu Choi\thanks{bchoi12@gsu.edu} }
\author[2]{Yoonsang Lee\thanks{yoonsang.lee@dartmouth.edu}}
\affil[1]{Department of Mathematics and Statistics, Georgia State University}
\affil[2]{Department of Mathematics, Dartmouth College}
\date{}

\begin{document}
		\maketitle
		
		\begin{abstract}
			In data assimilation, an ensemble provides a way to propagate \yl{the} probability density of a system described by a nonlinear prediction model. Although a large ensemble size is required for statistical accuracy, the ensemble size is typically limited to a small number due to the computational cost of running the prediction model, which leads to a sampling error. Several methods, such as localization and inflation, exist to mitigate the sampling error, often requiring problem-dependent fine-tuning and design. This work introduces a nonintrusive sampling error mitigation method that modifies the ensemble to ensure a smooth turbulent spectrum. It turns out that the ensemble modification to satisfy the smooth spectrum leads to inhomogeneous localization and inflation, which apply \yl{spatially} varying localization and inflation levels at different locations. The efficacy of the new idea is validated through a suite of stringent test regimes of the Lorenz 96 turbulent model.
		\end{abstract}
		
	\section{Introduction}
A wide range of problems in engineering and science involve estimating variables through a prediction model. As such prediction models are typically imperfect and involve uncertainty, observation data provides additional information to complement the insufficient information in the prediction model. Data assimilation utilizes Bayesian inference to provide statistically accurate and consistent estimation, combining uncertain predictions with observation data. The Kalman filter \cite{KF} provides an optimal estimation method under a linear prediction model with linear observation and Gaussian noises. In many applications, however, the prediction model is usually complicated and nonlinear, and thus, in general, no analytic update formula such as the Kalman filter exists for such nonlinear problems. 

Ensemble-based Kalman filters have attracted many researchers due to their nonintrusive approach to propagating and quantifying uncertainty related to a nonlinear prediction model. Ensemble Kalman filter (EnKF; \cite{EnKF}) uses an ensemble to sample probability distributions and their propagation along with the perturbed observation approach to assimilate observation data. The ensemble square root filters \cite{ESRF}, such as ensemble adjustment \cite{EAKF} or ensemble transform \cite{wang2003comparison} Kalman filters, use deterministic transformation (or adjustment) of an ensemble to satisfy the posterior statistics by the Kalman update. 

Despite a wide range of successful applications of ensemble-based Kalman filters, there are several sources of errors in ensemble-based Kalman filters: i) linear update (Gaussian \yl{assumption}), ii) prediction model error, iii) sampling error, etc. 
As in the standard Kalman filter, the ensemble-based Kalman filter uses a Gaussian assumption on the prior, which can lead to assimilation errors for a system with non-Gaussian statistics. There is a recent study on nonlinear update ensemble filtering \cite{youssef2022}, reducing the intrinsic bias of the ensemble Kalman filter using transport maps. Another class of data assimilation methods, particle filters \cite{PF,CPF}, use a nonparametric representation of a probability distribution through particles, enabling the natural treatment of non-Gaussian distributions. 

Also, when the prediction model is imperfect, the model's prior prediction will contain bias. When such bias is large, the prior distribution will deviate far from the unknown true distribution and thus degrade the overall data assimilation performance. Covariance inflation \cite{EAKF, inflation} increases the spread of the ensemble to pull back the ensemble close to the observation (and thus to the unknown true distribution). In addition to increasing the accuracy of addressing the model error, various inflation methods are known to prevent filter divergence \cite{inflation, preventCFD}, which is related to the stability of ensemble data assimilation methods. 

Additionally, ensemble-based Kalman filters suffer from sampling errors due to \yl{their} small ensemble size. The computational cost for a high-dimensional complex system (such as turbulent systems) is tremendously high, \yl{so} the ensemble size is limited to a relatively small size compared to the dimension of the variable of interest, which leads to sampling errors.
Localization \cite{localization} is a widely used method to handle sampling errors. 
When the ensemble size is small, it is known that spurious correlations in the prior covariance matrix generate artificial updates. Localization deflates the sample correlations based on certain criteria, such as physical distance between the observed and unobserved components \cite{gaspari_cohn1998} or off-line statistical tests \cite{sampErrCorrection2012}.

Although inflation and localization are crucial in stabilizing and mitigating sampling errors of ensemble-based Kalman filters, the optimal choice of such approaches requires an intensive tuning process. Such an intensive tuning process limits the optimal selection of inflation and localization levels to be applied uniformly regardless of the point of interest to apply inflation and localization. In the case of inflation, the same inflation level applies to each component of the state variable. In the case of localization, the same localization shape applies regardless of the observation location. There is a recent study that implements adaptive localization using convolutional neural networks \cite{anderson2023}. However, such a method requires significant training data to learn the localization function.

The goal of the current study is a nonintrusive sampling error mitigation method for turbulent systems. The Wiener-Khinchin theorem \cite{Khinchin} shows a one-to-one correspondence between the autocorrelation of a stationary stochastic process and its spectrum. That is, the autocorrelation of a stationary stochastic process determines its spectrum, and vice versa. Thus, estimating a correct spectrum is crucial in providing information about the interrelation between different locations. For fully developed turbulent systems, the Kolmogorov theory \cite{mitturbulence, chorin2013vorticity} shows that its spectrum follows a certain decaying rate, implying the smoothness of the spectrum as a function of the wavenumber. As this result is derived from statistical averaging, we expect the spectrum of the ensemble to be smooth under a large ensemble limit. Following this rationale, our proposed method modifies the ensemble to ensure its spectrum satisfies certain smoothness, which leads to sampling error mitigation.

The rest of this paper is organized as follows. Section \ref{sec:reviewKF} reviews standard sampling error mitigation methods in the context of the ensemble transform Kalman filter (ETKF). Ensemble modification based on the smooth spectrum constraint is described in detail in Section \ref{sec:ourmethod}, along with related theories. In Section \ref{sec:numerics}, we provide a suite of numerical test results using Lorenz 96 with various forcing constants ranging from weak turbulent to strong turbulent regimes, which supports the efficacy of the new method in mitigating sampling errors. We finish the study with discussions and future directions in Section \ref{sec:conclusion}.

\section{Sampling error in data assimilation}\label{sec:reviewKF}

This section reviews the ensemble-based Kalman filter along with problem setups of interest. Although spectrum smoothing can be applied to any ensemble-based Kalman filter, we will demonstrate the efficacy of the spectrum smoothing-based sampling error mitigation method in the context of, in particular, the ensemble transformation Kalman filter \cite{ETKF} in Section \ref{sec:numerics}. Thus, for the completeness of the current study, we focus on ETKF in reviewing ensemble-based Kalman filters. This section also discusses standard inflation and localization methods, \yl{which are} crucial in improving and stabilizing data assimilation methods. It turns out in Section \ref{sec:numerics} that the sampling error mitigation based on spectrum \rv{ smoothing} plays a role in \rv{ inhomogeneous} inflation and localization, which applies inflation and localization inhomogeneously based on each component of the system.

\subsection{Ensemble Transformation Kalman Filter} \label{subsec:problem}

Data assimilation incorporates observation data to improve the statistical accuracy of a dynamical system with uncertainty. In the current study, we assume the state variable $\bm{u}(t) \in \mathbb{R}^N$ follows a certain evolution (or prediction) map $\Psi:\mathbb{R}^N \rightarrow \mathbb{R}^N$ which we assume to be known. Additionally, we assume there are observation data $\bm{y}=h(\bm{u})\in\mathbb{R}^M$ available at discrete times $t_j=j\Delta t$, which are generally incomplete and noisy. For the simplicity of the argument, we assume that the time interval of the evolution operator matches with the observation time interval and thus use the subscript to represent the time index. Where $\epsilon_{j+1}$ represents the observation error, which is assumed to be unbiased Gaussian with a known covariance matrix $\Gamma$, we have the following discrete dynamical system and an observation model
\begin{eqnarray}
	\bm{u}_{j+1}&=&\Psi(\bm{u}_j),\\
	\bm{y}_{j+1}&=&h(\bm{u}_{j+1})+\epsilon_{j+1}.
\end{eqnarray}

When the system involves uncertainty (e.g., in the modeling of the dynamical systems or in the system parameters, such as an initial value), assimilating the observation data can improve statistical accuracy in estimating the state variable $\bm{u}$. Under linear assumptions on both the evolution map $\Psi(\bm{u}_j)=M\bm{u}_j$ and the observation map $h(\bm{u}_j)=H\bm{u}_j$, the Kalman filter \cite{KF} update provides the optimal statistical estimation for the mean $\bm{m}_j$ and covariance $C_{j}$, which involves two steps: i) prediction, and ii) assimilation as follows

\textbf{i) Prediction:}
\begin{eqnarray}
	\bm{\hat{m}}_{j+1}&=&M\bm{m}_j^{posterior} \\
	\hat{C}_{j+1}&=&MC_{j}M^T\\
\end{eqnarray}

\textbf{ii) Assimilation:}
\begin{eqnarray}
	\bm{m}_{j+1}&=&\bm{\hat{m}}_{j+1}+K_{j+1}\left(\bm{y}_{j+1}-H\bm{\hat{m}}_{j+1}\right)\\
	C_{j+1}&=&(I-K_{j+1}H)\hat{C}_{j+1}
\end{eqnarray}
where $K_{j+1}$ is the Kalman gain
\begin{equation}
	K_{j+1}=\hat{C}_{j+1}H^T(H\hat{C}_{j+1}H^T+\Gamma)^{-1}.
\end{equation}

When the evolution map $\Psi$ is nonlinear, typical in a wide range of applications in science and engineering, including numerical weather prediction, the prediction and its corresponding propagation of uncertainty becomes a nontrivial computational problem. Ensemble-based Kalman filters, such as ensemble Kalman filter \cite{EnKF}, ensemble transformation \cite{ETKF}, and ensemble adjustment \cite{EAKF} Kalman filters, use an ensemble (or a sample of the system) of size $K$, $\{\bm{u}^{k}_{j}\}_{k=1}^{K}$, to propagate the system
\begin{equation}
	\bm{u}_{j+1}^{(k)}=\Psi(\bm{u}_{j}^{(k)}), \quad k=1,2,...,K.
\end{equation}
Once the ensemble has evolved, they use ensemble statistics for the prior mean and covariance
\begin{eqnarray}
	\hat{\bm{m}}_{j+1} &=& \frac{1}{K} \sum_{k=1}^K  \rv{ \hat{\bm{u}}_{j+1}^{(k)}}\\
	\hat{C}_{j+1} &=& \frac{1}{K-1} \sum_{k=1}^K ( \rv{ \hat{\bm{u}}_{j+1}^{(k)}} - \hat{\bm{m}}_{j+1})
	( \rv{ \hat{\bm{u}}_{j+1}^{(k)}}  - \hat{\bm{m}}_{j+1})^T.
\end{eqnarray}	

In the assimilation step, the mean is updated as in the standard Kalman update using the ensemble mean and covariance
\begin{eqnarray}
	\bm{m}_{j+1}&=&\bm{\hat{m}}_{j+1}+K_{j+1}\left(\bm{y}_{j+1}-H\bm{\hat{m}}_{j+1}\right) \label{eq:post_mean}\\
	K_{j+1}&=&\hat{C}_{j+1}H^T(H\hat{C}_{j+1}H^T+\Gamma)^{-1}.  \label{eq:kalmanGain}  
\end{eqnarray}
The aforementioned methods differ based on how they modify the ensemble to satisfy the posterior statistics. In the current study, we focus on the ensemble transformation Kalman filter, which uses a transformation matrix to satisfy the posterior covariance matrix of the Kalman update. 

In particular, \yl{following the implementation of ETKF in \cite{Law15}, in order to transform the perturbation part}, $\hat{X}_{j+1}$, of the ensemble defined as
\begin{equation}
	\hat{X}_{j+1} = \frac{1}{\sqrt{K-1}}
	\begin{bmatrix} \hat{\bm{u}}_{j+1}^{(1)} - \hat{\bm{m}}_{j+1}, \cdots,
		\hat{\bm{u}}_{j+1}^{(K)} - \hat{\bm{m}}_{j+1}\end{bmatrix},
\end{equation}
we calculate
\begin{equation}
	T_{j+1}= [I + (H \hat{X}_{j+1})^T \Gamma^{-1} (H \hat{X}_{j+1})]^{-1},  \label{eq:tmat}
\end{equation}
\rv{ which is applied for updating the perturbation part $\hat{X}_{j+1}$. The symmetric square root of $T_{j+1}$ is multiplied to $\hat{X}_{j+1}$ to \rv{obtain}
\begin{equation}
	X_{j+1}= \hat{X}_{j+1} T_{j+1}^{1/2}  \label{eq:postmat}
\end{equation}
in {\bf the assimilation step} while the ensemble mean, $\bm{m}_{j+1}$, is updated separately as in \eqref{eq:post_mean}.}
By adding the posterior mean $\bm{m}_{j+1}$ to each column of $X_{j+1}$, \rv{ we obtain the posterior ensemble, $\bm{u}_j^{(k)} = \bm{m}_{j+1} +  \bm{x}_{j+1}^{(k)}$, where $\bm{x}_{j+1}^{(k)}$ is the $k$-th column of $X_{j+1}$.}

In addition, for the evolution map $\Psi$ to be nonlinear, its computational cost to solve is typically high and thus requires tremendously large computing resources. Such high cost in evolving each ensemble member limits the ensemble size, typically much less than the state dimension, $K\ll N$, leading to various sampling errors.

\subsection{Localization and inflation} \label{sec:err_correction}
When the ensemble size is not sufficiently large, there are several issues that degrade the performance of the ensemble-based Kalman filters. The correlation, which is crucial in updating unobserved components using data, can have a significantly large value between uncorrelated variables, called `spurious correlation' \cite{localization}. Also, when the ensemble size is smaller than the state dimension ($K<N$), the sample covariance matrix suffers from rank deficiency. In the case of ETKF, for example, rank deficiency can result in failure to obtain the transformation matrix. To address such issues, various localization \cite{localization} methods deflate the correlation information to suppress spurious correlation values. 

A distance-dependent localization factor has been designed for a physical system, including the widely used localization factor by Gaspari and Cohn \cite{gaspari_cohn1998}, which is a fifth-order rational function with compact support. The rationale behind distance-dependent localization is that two variables sufficiently far apart must be weakly correlated while adjacent variables are strongly correlated. 
Given a distance-based localization factor function $g(\cdot):\mathbb{R}\to [0,1]$, localization calculates a matrix $L$ whose component $i$th row and $j$th column is given by $g(\|i-j\|_p)$ where $\|i-j\|_p$ is a distance between the $i$-th and $j$-th component of the state variable. In the one-dimensional physical space, for example, each row of $L$ is translation invariant. That is, the localization factor is homogeneous across different rows. Once the localization matrix $L$ is designed, the covariance matrix is modified as follows
\begin{equation}
	\hat{C}_{j} \leftarrow L \odot \hat{C}_{j},
\end{equation}
where $\odot$ represents a Hadamard product.

Localization performance is sensitive to several parameters, such as the domain of influence. In the Gaspari-Cohn localization, the halfwidth (often denoted as $c$) determines the size of the localization support. Due to the sensitive tuning process in determining the halfwidth parameter, it is a challenging task to design different localization factors in each row of $L$. There are several other approaches aimed at mitigating sampling error. The sampling error correction in \cite{sampErrCorrection2012} produces a lookup table of localization factors based on offline samples of various correlation values. Also, recent work includes convolutional neural networks that enable adaptive localization \cite{anderson2023}, which requires training data to learn the localization function.

In addition to localization, inflation is crucial in stabilizing and improving the accuracy of ensemble-based Kalman filters \cite{inflation, adaptiveinflation, preventCFD}. Inflation adds uncertainty in the prior ensemble to compensate for the prior bias from a model error or sampling error. Additive inflation adds noise to each ensemble member, while multiplicative inflation spreads the ensemble to cover a wider range of prior values. In other words, both additive and multiplicative inflations modify the prior ensemble so that its corresponding covariance has a larger spread after the modification
\begin{eqnarray}
	\mbox{additive: } \hat{C}_{j}&\leftarrow& \hat{C}_{j}+ \alpha I, \quad \alpha\geq0\label{eq:addinf}\\
	\mbox{multiplicative: } \hat{C}_{j}&\leftarrow&\rho\hat{C}_{j}, \quad \rho\geq1\label{eq:mulinf}.
\end{eqnarray}
As in the localization, optimal inflation level selection remains a tuning process. Due to the complexity of choosing the optimal inflation level, it becomes challenging to choose various inflation levels based on different components. Several adaptive inflation methods change the inflation level in the spatiotemporal domain \cite{inflation, adaptiveinflation}, which also requires sufficient observation data statistics to design the inflation. Another adaptive inflation method uses the innovation statistics \cite{preventCFD}. Although this method improves in various stringent test cases, its performance can degrade when the observation dimension is significantly smaller than the state dimension.

\section{Sampling error mitigation based on smoothness constraint} \label{sec:ourmethod}

In ensemble-based data assimilation, statistics of ensemble and observations determine the evolution of probability distribution. Therefore, the ensemble size (i.e., sample size) and the estimation of correlation between \yl{the} state variable's components affect the accuracy of probability estimation. \yl{We} pay attention to Kolmogorov's theory \cite{chorin2013vorticity}  which provides the statistical characterization of turbulence, which yields the constraints that a large ensemble statistics would satisfy. Such constraint is observed in the \yl{turbulence spectrum}, which decays at a certain rate and hence is smooth. Moreover, the spectrum has the one-to-one correspondence with the autocorrelation by  Khinchin Theorem \cite{chorin2013vorticity}. Accordingly, the accuracy of \yl{the} correlation estimate is expected to be improved by modifying the ensemble members in a way that its spectrum becomes closer to what the spectrum of a larger ensemble is supposed to be. Consequently, statistics of the updated ensemble better describe the probability distribution of a state variable without increasing the ensemble size.

Kolmogorov's theory and Khinchin Theorem are reviewed in Section \ref{sec:reviewTurbulence}. We also \yl{introduce} how we define the spectrum of an ensemble. For simplicity, we use a one-dimensional turbulence system in this section. However, \yl{we want to emphasize} that this definition can be extended to high dimensional ones, which is skipped in this paper. To support our approach using the smooth constraint of a spectrum, we compare the spectra of small and large ensembles, and also compare those of small ensembles with and without our spectrum smoothing. The results of this numerical comparison and the algorithmic description of the spectrum smoothing are introduced in Section \ref{subsec:samp_err_mitigation}.

\subsection{Turbulence Theory} \label{sec:reviewTurbulence}

A fully developed turbulent flow is very sensitive to a small perturbation and shows chaotic and random behavior. Turbulence \rv{theory} describes such turbulent flow using the statistical theory, which reveals more robust properties of turbulence. Despite its random behavior, a turbulent flow consists of large- to small-scale turbulent motions which interact and transfer the energy to each other with certain statistical pattern (i.e. \rv{homogeneity}). Such statistical pattern is well-described in Kolmogorov theory \cite{chorin2013vorticity}. This energy transfer is called ``energy cascade" whose name comes from the fact that the corresponding \rv{energy} spectrum decays with a \rv{certain} rate.

We briefly summarize the turbulence theory to introduce Kolmogorov theory and Khinchin Theorem. We start by introducing the energy spectrum using a one-dimensional random flow. Its extension to the high-dimensional flow is available in \cite{chorin2013vorticity}. 
Let $u(x)$ be a random flow field, and under the assumption of periodicity, we can write $u$ as a Fourier expansion as follows:
\begin{equation}
	u(x) = \sum c_{\omega} e^{i {\omega} x},
\end{equation}
with a wavenumber $\omega$.
We define the spectral distribution, $F$, of $u$ as:
\begin{equation}
	F(\bar{\omega}) = \sum\limits_{|\omega| < {\bar{\omega}}} \langle |c_{\omega}|^2 \rangle,
\end{equation}
where $\langle \eta \rangle$ is defined to be the expected value of $\eta$, i.e.,  $\langle \eta \rangle :=\int_{\Omega}\eta dP$ with a probability space, $(\Omega, \mathcal{B}, P)$, and $|\cdot|$ is a Euclidean distance. Under the assumption that $F$ is differentiable,  we obtain
\begin{equation}
	F'({\bar{\omega}}) = \phi ({\bar{\omega}}),
\end{equation}
where $\phi$ is the spectral density of $u$. Finally, we define the energy spectrum, $E(\omega)$, of $u$ as the following:
\begin{equation}
	E(\omega) = \frac{1}{2} \int_{|\tau|=\omega} \phi (\tau) d \tau.
\end{equation}
According to the Kolmogorov theory, the energy spectrum, $E(\omega)$, of homogeneous flow follows:
\begin{equation}
	E(\omega) = C \epsilon^{2/3} \omega^{-5/3},
\end{equation}
where $C$ is a dimensionless constant, and $\epsilon = \frac{d}{dt}\langle u^2 \rangle$ is the rate of energy dissipation. This indicates that the energy spectrum of a random flow field decays with a certain rate as $\omega$ increases so that it is smooth. 

On the other hand, Wiener-Khinchin Theorem or  Khinchin Theorem relates the energy spectrum and the autocorrelation: one-to-one correspondence of each other. Khinchin Theorem states that the autocorrelation function of a stationary stochastic process has a spectral decomposition given by the power spectral density of the process \cite{chatfield2004timeseries}.
Instead of considering the autocorrelation, however, we can apply the theorem in terms of  the correlation in the physical space
if a turbulent flow is homogeneous and isotropic.
The correlation function of $u(x)$ is
\begin{equation}
	R(x_1,  x_2)  = \langle (u(x_1) - m(x_1))(u(x_2) - m(x_2)) \rangle
\end{equation}
where $m(x)$ is a mean of $u(x)$. Then, Khinchin \rv{Theorem} states that
\begin{theorem*}[Khinchin Theorem]
	For the correlation function $R(x)$ of a random flow field, $u(x)$, which has translation invariant mean and correlation function and also satisfies the condition 
	\begin{equation}
		\langle |u(x+h) - u(x)|^2 \rangle \rightarrow 0 ~ as ~ h \rightarrow 0,
	\end{equation}
	it is necessary and sufficient that it has a representation of the form 
	\begin{equation}
		R(x) = \int_{-\infty}^{\infty} \exp(ix\omega) dF(\omega),
	\end{equation}
	where $F(\omega)$ is a spectral distribution of $u(x)$. 
\end{theorem*}
Khinchin Theorem implies that for a given $F(\omega)$, $R(x)$ is the inverse Fourier transform of the spectral density function, $\phi(\omega)= F'(\omega)$, and therefore, the correlation function has a one-to-one correspondence with the spectral density function and the energy spectrum as well. 

Based on the implication of Khinchin Theorem, we expect that the better estimation of energy spectrum would improve the estimation of the correlation, which is crucial in the data assimilation.
From this rationale, \yl{we} suggest improving the energy spectrum estimation by adjusting the mean spectrum of ensemble. Accordingly, we expect that this adjustment increases the accuracy of correlation/covariance estimations. We estimate the energy spectrum, $E(\omega)$, by the mean power spectrum of the ensemble which is defined as follows: 
\begin{equation}
	\frac{1}{K} \sum_{k=1}^K |\mathcal{F} ( \bm{u}^{(k)} )[\omega]|^2,
\end{equation}
where $\bm{u}^{(k)}$ is the $k$th ensemble member, and $\mathcal{F}(\bm{u})[\omega]$ represents the discrete Fourier transform of $\bm{u}$ at a wavenumber $\omega$.
The new sampling error mitigation using the spectrum smoothing is introduced in detail in Section \ref{subsec:samp_err_mitigation}.

\subsection{Spectrum smoothing} \label{subsec:samp_err_mitigation}

Inspired by the smooth characteristics of the energy spectrum, we adjust \yl{the} ensemble so that its mean power spectrum becomes closer to the energy spectrum without increasing the ensemble size. When the ensemble size is small, the mean power spectrum is typically nonsmooth due to the sampling error. Therefore, we obtain a smoother spectrum by taking the convolution of a mean power spectrum with a smoothing kernel. This smoother spectrum is expected to be closer to the energy spectrum, and accordingly, we adjust the ensemble so that its mean spectrum \yl{becomes} the smoother spectrum. 
\yl{Consequently}, the adjusted ensemble has \yl{a} more accurate mean spectrum, and therefore, it is expected to reflect the underlying probability distribution more accurately.

To elaborate on the spectrum smoothing, we define a discrete convolution of a given vector $\phi$ with a smoothing kernel, $\kappa_{\sigma}$, as follows:
\begin{equation}
	(\phi * \kappa_{\sigma})[\omega] = \sum_{\tau= -\infty}^{\infty} \phi [\omega - \tau]\kappa_{\sigma}[\tau],
	\label{eq:smoothing}
\end{equation}
where $\sigma$ is a kernel width scaling factor.
The ensemble is updated by rescaling each ensemble member so that the updated ensemble's mean spectrum becomes the smoothed spectrum. Specifically,
we update the prior ensemble  $\hat{\bm{u}}_{j+1}^{(k)} = \hat{\bm{m}}_{j+1} + \hat{\bm{x}}_{j+1}^{(k)}$,	where $\hat{\bm{m}}_{j+1} = \frac{1}{K} \sum_{k=1}^K \hat{\bm{u}}_{j+1}^{(k)}$ is the mean ensemble by rescaling the perturbation, $\hat{\bm{x}}_{j+1}^{(k)}$, in the Fourier domain. By doing so, the ensemble mean remains the same, but the perturbations are modified to improve the energy spectrum estimate. 
Consequently, the updated prior ensemble's discrete Fourier transform becomes
\begin{equation}
	\mathcal{F}(\hat{\bm{m}}_{j+1})[\omega]  +  \alpha[\omega] \mathcal{F}(\hat{\bm{x}}_{j+1}^{(k)})[\omega]
\end{equation}
with a multiplicative rescaling factor $\alpha[\omega]$.
This rescaling factor $\alpha[\omega]$ for each wavenumber $\omega$ is designed as follows:
\begin{equation}	
	\alpha[\omega] = \sqrt{ \frac{S_{\sigma} \left( \frac{1}{K} \sum_{k=1}^K |\mathcal{F} ( \hat{\bm{u}}_{j+1}^{(k)} )|^2 \right)[\omega]  -  |\mathcal{F}(\hat{\bm{m}}_{j+1})[\omega]|^2}{ \frac{1}{K} \sum_{k=1}^K  |\mathcal{F}(\hat{\bm{x}}_{j+1}^{(k)})[\omega]|^2} }.
	\label{eq:alpha_eq2}
\end{equation}
$\alpha[\omega]$ in \eqref{eq:alpha_eq2} is derived \yl{from making} the updated mean power spectrum in \eqref{eq:alpha_eq} below equal to the smoothed mean power spectrum, i.e., the second equality below holds as in the following:
\begin{align}
	S_{\sigma} \left( \frac{1}{K} \sum_{k=1}^K |\mathcal{F} ( \hat{\bm{u}}_{j+1}^{(k)} )|^2 \right)[\omega]
	&:= \begin{cases} \left( \left( \frac{1}{K} \sum_{k=1}^K |\mathcal{F} ( \hat{\bm{u}}_{j+1}^{(k)} )|^2 \right)  * \kappa_{\sigma} \right) [\omega]
		\\ \hspace{2em}{\rm ~if~} \left( \left( \frac{1}{K} \sum_{k=1}^K |\mathcal{F} ( \hat{\bm{u}}_{j+1}^{(k)} )|^2 \right)  * \kappa_{\sigma} \right)[\omega] \geq |\mathcal{F}(\hat{\bm{m}}_{j+1})[\omega]|^2   \\ 
		|\mathcal{F}(\hat{\bm{m}}_{j+1})[\omega]|^2 {\rm ~otherwise}
	\end{cases} \label{eq:smoothingOp}\\
	&= \frac{1}{K} \sum_{k=1}^K |\mathcal{F}(\hat{\bm{m}}_{j+1})[\omega]  +  \alpha[\omega] \mathcal{F}(\hat{\bm{x}}_{j+1}^{(k)})[\omega]|^2,
	\label{eq:alpha_eq}
\end{align}
$S_{\sigma}(\cdot)$ is a smoothing operator with a kernel width scaling factor $\sigma$. The operator is equivalent to the convolution with a smoothing kernel at a wavenumber $\omega$ if the smoothed spectrum is no less than the mean spectrum. If it is less, we set the smoothing operator to equal the mean spectrum for such \yl{a} wavenumber. The smoothing operator $S_{\sigma}(\cdot)$ is defined in this way as in \eqref{eq:smoothingOp}  to make the rescaling factor $\alpha[\omega]$ to be a real function at every $\omega$. Otherwise, $\alpha[\omega]$ in \eqref{eq:alpha_eq2} cannot be well-defined. Using the fact that the sum of perturbations, $ \hat{\bm{x}}_{j+1}^{(k)}$, is zero, and accordingly, $\sum_{k=1}^K \mathcal{F} ( \hat{\bm{x}}_{j+1}^{(k)} ) = 0$.
Then, the updated ensemble mean power spectrum, \eqref{eq:alpha_eq}, is simplified as follows:
\begin{equation}
	S_{\sigma} \left( \frac{1}{K} \sum_{k=1}^K |\mathcal{F} ( \hat{\bm{u}}_{j+1}^{(k)} )|^2 \right) [\omega]
	= |\mathcal{F}(\hat{\bm{m}}_{j+1})[\omega]|^2  +  \frac{\left( \alpha[\omega] \right)^2}{K} \sum_{k=1}^K  |\mathcal{F}(\hat{\bm{x}}_{j+1}^{(k)})[\omega]|^2.\\
	\label{eq:rescale}
\end{equation}	
By solving  \eqref{eq:rescale} for $\left( \alpha[\omega] \right)^2 $, we obtain:	
\begin{equation}
	\left( \alpha[\omega] \right)^2 = \frac{S_{\sigma} \left( \frac{1}{K} \sum_{k=1}^K |\mathcal{F} ( \hat{\bm{u}}_{j+1}^{(k)} )|^2 \right)[\omega] -  |\mathcal{F}(\hat{\bm{m}}_{j+1})[\omega]|^2}{ \frac{1}{K} \sum_{k=1}^K  |\mathcal{F}(\hat{\bm{x}}_{j+1}^{(k)})[\omega]|^2}. \\
\end{equation}	
We consistently choose the nonnegative sign for determining $\alpha[\omega]$ to avoid switching the phase, and therefore, we obtain the rescaling function $\alpha[\omega]$ in \eqref{eq:alpha_eq2}.
Here, the smoothing operator $S_{\sigma}$ is designed to ensure that the numerator in \eqref{eq:alpha_eq2} becomes nonnegative. Therefore, the smoothed power spectrum is no less than the mean power spectrum of the prior ensemble as in \eqref{eq:smoothingOp}. We emphasize that \eqref{eq:alpha_eq2} is computed component-wisely for each wavenumber $\omega$. 
Eventually, the updated ensemble can be obtained by taking the inverse discrete Fourier transform of each ensemble member rescaled in the Fourier domain, i.e., 
\begin{equation}
	\mathcal{F}^{-1}\left(   \mathcal{F}(\hat{\bm{m}}_{j+1})[\omega]  +  \alpha[\omega] \mathcal{F}(\hat{\bm{x}}_{j+1}^{(k)})[\omega] \right).
\end{equation}

\begin{figure}[h]
	\centering
	\subfloat[Mean power spectrum at $t=48.8$ \\ {\it without} smoothing]
	{\includegraphics[width=.45\textwidth]{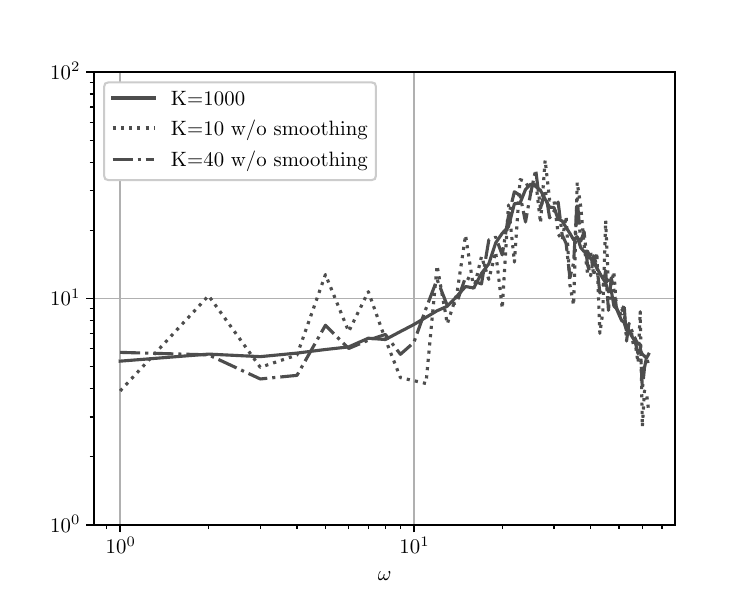}\label{fig:l96sampleSpec0}}
	\subfloat[Mean power spectrum at $t=48.8$ \\ {\it with} smoothing]
	{\includegraphics[width=.45\textwidth]{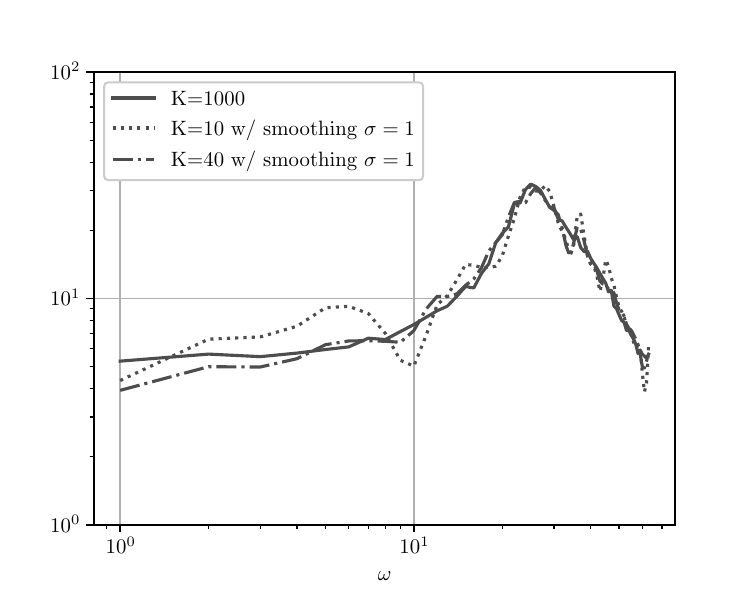}\label{fig:l96sampleSpec2}}
	\caption{Lorenz 96 mean power spectrum comparison while varying the ensemble size: \\ (a) spectrum smoothing is {\it not} applied and (b) spectrum smoothing is applied}
	\label{fig:l96sampleSpec}
\end{figure}

After updating the ensemble, the mean power spectrum is smoother than before, and the resulting smoothed spectrum would be more similar to the spectrum of the ensemble in \yl{a} larger size. In this way, we expect the updated ensemble to \yl{better} represent the underlying probability distribution of the state variables. Since the data assimilation tends to be more effective when the ensemble size is larger, we expect that the sampling errors in the new ensemble whose spectrum is closer to the larger ensemble's spectrum are mitigated by the spectrum smoothing. We demonstrate the suggested spectrum smoothing \yl{effect} by applying it to the ensemble for \yl{the} Lorenz 96  model with a forcing constant $F=8$.  The model is solved by the 4th-order Runge-Kutta method, and the ensemble is propagated by the free-run simulation, i.e., the assimilation using new observations is skipped. We observe in Fig \ref{fig:l96sampleSpec} that the mean power spectrum tends to be smoother as the ensemble size increases.
Furthermore, when applying the spectrum smoothing, the mean power spectrum becomes closer to the one from the ensemble in size 1000. This demonstrates that the ensemble update using the spectrum smoothing works as intended. In Section \ref{sec:numerics}, the performance of ETKF with the updated ensemble with a smoother spectrum will be tested to validate the sampling mitigation effect.

Algorithm \ref{alg:ETKFsmoothing} describes the ETKF with the spectrum smoothing. The parameters, $c$, $\rho$, and $\sigma$, are chosen for localization, inflation, and spectrum smoothing level, respectively. ETKF has two iterative steps: prediction and assimilation. The prediction step is the same as the classic ETKF, while the spectrum smoothing is performed repeatedly before each assimilation step. Right before the assimilation step, the {\it prior} ensemble $\hat{\bm{u}}_{j+1}$ is updated by rescaling its \rv{purturbation} part using the rescaling factor $\alpha$. After the ensemble update, localization and inflation are still necessary since they alleviate other sources of errors such as model errors, observation errors, rank deficiency, etc., as discussed in Section  \ref{sec:reviewKF}. \rv{In the algorithm, multiplicative inflation is applied by multiplying the inflation factor $\sqrt{\rho} \geq 1$ to the ensemble perturbation. That is, the variance is increased by a uniform factor $\rho$. To suppress the spurious correlation, on the other hand, covariance localization is applied by taking the Hadamard product of the prior covariance matrix and the localization matrix, $L$. Each element $L_{ij}$ of $L$ is determined by the Gaspari-Cohn function value based on the spatial distance, $\|i-j\|$,  and the parameter $c$ scaling the support of \yl{the} Gaspari-Cohn function. As the distance is farther, a smaller $L_{ij}$ \yl{multiplies} to the element of the covariance matrix. }
Therefore, in Section \ref{sec:numerics}, the performance of Algorithm \ref{alg:ETKFsmoothing} is tested along with the classic ETKF with localization and inflation for a fair comparison. Though the smoothing spectrum idea is implemented within ETKF in this work, it can be used in any ensemble-based method to update the ensemble to have a smoother spectrum.

\begin{algorithm}[!t]
	\caption{ETKF with the spectrum smoothing}\label{alg:ETKFsmoothing}
	\begin{algorithmic}
		\State Initialize the ensemble $\bm{u}_{j}^{(k)}$. Choose parameters $c$, $\rho$, and $\sigma$ 
		\For{$j = 1$ to $J-1$}	
		\State ({\it Prediction}) compute the prior ensemble $\hat{\bm{u}}_{j+1}^{(k)} = \Psi(\bm{u}_{j}^{(k)})$ 
		for $k = 1, \cdots, K$.
		\If{observations are available at $j+1$-th time step}
		\State ({\it Assimilation})
		\State Update the prior ensemble using the spectrum smoothing with $\alpha$ in \eqref{eq:alpha_eq2}:
		\[\hat{\bm{u}}_{j+1}^{(k)} \leftarrow \mathcal{F}^{-1}(\mathcal{F}(\hat{\bm{m}}_{j+1})  +  \alpha \mathcal{F}(\hat{\bm{x}}_{j+1}^{(k)}))\]
		\State Inflate the prior ensemble spread: 
		$\hat{\bm{u}}_{j+1}^{(k)} \leftarrow \hat{\bm{m}}_{j+1} 
		+ \sqrt{\rho} \hat{\bm{x}}_{j+1}^{(k)}$	 
		\State Update the prior covariance matrix with a localization: $\hat{C}_{j+1} \leftarrow L_c \odot \hat{C}_{j+1}$.	
		\rv{
		\State Compute the Kalman gain $K_{j+1}$ and transform matrix $T_{j+1}$ as in \eqref{eq:kalmanGain} and \eqref{eq:tmat}.
		\State Obtain the posterior ensemble mean, $\bm{m}_{j+1}$, and perturbation,  $X_{j+1}$ , as in  \eqref{eq:post_mean}  and \eqref{eq:postmat}.}
		\EndIf
		\EndFor	
	\end{algorithmic}
\end{algorithm}



In the suggested spectrum smoothing, we smooth out the mean power spectrum of the entire ensemble and rescaling the perturbation part of the ensemble in the Fourier domain. We also tried another approach for the comparison: smoothing the spectrum of the entire ensemble as before but rescaling the entire ensemble instead of the perturbation part, i.e., designing another rescaling factor $\tilde{\alpha}[\omega]$ so as to satisfy the following:
\begin{equation}
	S_{\sigma} \left( \frac{1}{K} \sum_{k=1}^K |\mathcal{F} ( \hat{\bm{u}}_{j+1}^{(k)} )|^2 \right)[\omega] = \frac{1}{K} \sum_{k=1}^K | \tilde{\alpha}[\omega] \mathcal{F}(\hat{\bm{u}}_{j+1}^{(k)})[\omega]|^2,
\end{equation}
where $\tilde{\alpha}[\omega]$ is multiplied to each ensemble member in the Fourier domain, $\mathcal{F}(\hat{\bm{u}}_{j+1}^{(k)})[\omega]$, instead of the perturbation part of the ensemble in the Fourier domain, $\mathcal{F}(\hat{\bm{x}}_{j+1}^{(k)})[\omega]$. $\tilde{\alpha}[\omega]$ can be explicitly obtained as the following
\begin{equation}	
	\tilde{\alpha}[\omega] = \sqrt{ \frac{S_{\sigma} \left( \frac{1}{K} \sum_{k=1}^K |\mathcal{F} ( \hat{\bm{u}}_{j+1}^{(k)} )|^2 \right)[\omega]}{|\mathcal{F}(\hat{\bm{m}}_{j+1})[\omega]|^2 + \frac{1}{K} \sum_{k=1}^K  |\mathcal{F}(\hat{\bm{x}}_{j+1}^{(k)})[\omega]|^2} }.
	\label{eq:tilde_alpha_eq}
\end{equation}
This approach changes the mean ensemble, $\hat{\bm{m}}_{j+1}$, after the spectrum smoothing since it rescales the entire ensemble.
Consequently, the former suggested spectrum smoothing shows more stable prediction \rv{accuracy} than the later method in comparison when it is applied to ensemble data assimilation.

In the suggested spectrum smoothing, the perturbation part of the ensemble is rescaled in the Fourier domain, resulting in the convolution in the physical domain. This leads to inhomogeneous localization and inflation effect in addition to the sampling error mitigation, which will be demonstrated in Section \ref{sec:numerics}. Moreover, the rationale on the spectrum smoothing is based on the turbulence theory which applies to the fully developed turbulent flows. Therefore, the effectiveness of the spectrum smoothing differs depending on the turbulence strength, which will also be validated in Section \ref{sec:numerics} as well.

\section{Numerical experiments} \label{sec:numerics}

This section provides a suite of tests to validate the spectrum smoothing in mitigating the sampling errors. We use ETKF with inflation and localization as a baseline method to compare. Localization used the Gaspari-Cohn(GC) localization with various halfwidth $c$, ranging from 1 to 15. For inflation, we use multiplicative inflation with $\rho$ ranging from 1 to 1.2. Localization and inflation have been hand-tuned with an increment of 1 and 0.01, respectively. In measuring the overall performance, we use the time-averaged RMSE using the last 350 assimilation cycle results out of 1333 assimilation cycles. As the observation time interval $\Delta t_{obs}$ is 0.15, \yl{1333 assimilation cycles are} equivalent to running the whole system until $t = 200$.
For the spectrum smoothing, we find an optimal kernel width parameter $\sigma$ testing values from 0.1 to 1 with the increment of 0.1. Here, \yl{a larger $\sigma$ value implies stronger smoothing while a smaller $\sigma$ is weaker smoothing}.
\yl{As localization and inflation play a vital role in stabilizing filters and resolving rank deficiency even after mitigating the sampling error, we apply both localization and inflation after the spectrum smoothing.}

\subsection{Test model and test regimes}

As the test problem, we consider $N$-dimensional Lorenz-96 model \cite{lorenz96}
\begin{equation}
	\frac{du_n}{dt} = (u_{n+1} - u_{n-2}) u_{n-1} - u_n + F, ~n = 1,2, \cdots, N.
\end{equation}
The system is periodic, that is, $u_0 = u_N$ and $u_{-1} = u_{N-1}$, and we vary the forcing term $F$ to reach various turbulent regimes, i) weak turbulence ($F = 4$), ii) moderate turbulence ($F = 8$), and iii) strong turbulence ($F = 16$). The ODE system is solved by a 4-th order Runge-Kutta with a time step $\Delta t = 0.01$. As the observation interval $\Delta t_{obs}$ is 0.15, each prediction involves 15 iterations of evolving the system with the time step $\Delta t = 0.01$. For each test regime’s initial value, we used a constant value equal to $F$ and made a small perturbation in the first component. The actual initial value is obtained after running the regime \rv{sufficiently} long so that the model develops its corresponding turbulent characteristics.

\begin{figure}[!hbtp]
	\subfloat[F=4]
	{\includegraphics[width=.33\textwidth]{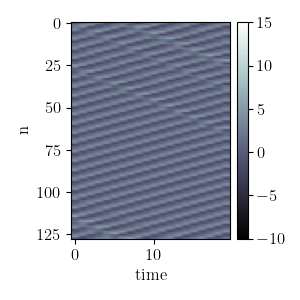}\label{fig:l96priorSol0_F4}}	
	\subfloat[F=8]
	{\includegraphics[width=.33\textwidth]{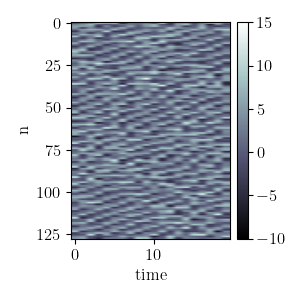}\label{fig:l96priorSol0}}
	\subfloat[F=16]
	{\includegraphics[width=.33\textwidth]{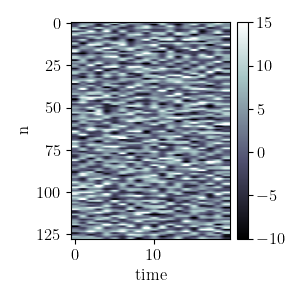}\label{fig:l96priorSol0_F16}}	
	\caption{Time series of 128-dimensional Lorenz 96 with varying $F \in\{4, 8, 16\}$}
	\label{fig:l96priorSol0_varyF}
\end{figure}

\begin{figure}[!hbtp]
	\centering
	{\includegraphics[width=.6\textwidth]{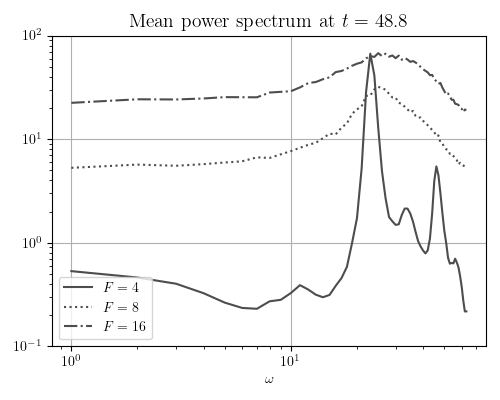}\label{fig:sampleSpecF}}
	\caption{Ensemble spectrum of 128-dimensional Lorenz 96 for various turbulent regimes. Each spectrum is obtained by using an ensemble of size $K=1000$ mean power spectrum comparison while varying $F$, each mean power spectrum is obtained by taking the mean from ensemble of size $K=1000$}
	\label{fig:sampleSpecCompF}
\end{figure}

\yl{The dimension of the system, $N$, is set to $128$ in all tests in order to have a fully developed turbulent spectrum, in particular, to have the decaying regime in the spectrum\footnote{\yl{We, in fact, applied the spectrum smoothing to the 40-dimensional case (results are not reported here), which is one of the standard test dimension. However, the performance improvement by the spectrum smoothing remains marginal as the corresponding spectrum of the 40-dimensional case has no significant smoothly decaying part.}}.} 
Figure \ref{fig:l96priorSol0_varyF} shows the free run (i.e., evolving the model without assimilating observation) time series of each test regime ($F=4$ (top), 8 (middle), and 16 (bottom)). To check the spectrum of each test regime, we add random perturbations to generate $K = 1000$ ensemble members and calculate the spectrum using the ensemble, which is shown in Figure \ref{fig:sampleSpecCompF}. As the model becomes more turbulent as $F$ increases, we find that the spectrum becomes smoother while its magnitude also increases (that is, the energy is increased too).

From these free runs, we find that the climatological standard deviation of each test regime is 1.854, 3.640, and 6.298. Using these values, we add observation noise in the data assimilation experiment whose standard deviation corresponds to \yl{10\% of the standard deviations}. Also, for the spatial resolution of the observation, we test i) 25\%, ii) 33\%, iii) 50\%, and iv) 100\% observations, which correspond to observing every i) three, ii) two, iii) other points while 100\% provides the full observation.

\rv{In order to measure the errors, the root mean square errors (RMSE) of classic ETKF and ETKF with the spectrum smoothing are obtained. RMSE at the $j$th time step is defined as the following,
\begin{equation}
	{\rm RMSE}_j = \sqrt{\frac{\| \bm{m}_{j} -\bm{u}_{j}\|_2^2}{N}},
\end{equation}
where $\bm{m}_{j}$ is the {\it posterior} ensemble mean, and $\bm{u}_{j}$ is the true solution at the $j$th time step. For comparing the performance, the average ${\rm RMSE}$ is computed by averaging ${\rm RMSE}_j$ over the last $350$ assimilation cycles out of 1333 cycles. Thus, the time-averaged {\it posterior} RMSE is:
\begin{equation}
	{\rm RMSE} = \frac{1}{350} \sum_{j=984}^{1333} {\rm RMSE}_{j}.
\end{equation}
}

\subsection{Test results}

In Section \ref{sec:inhomog}, we demonstrate that spectrum smoothing makes the correlation localized inhomogeneously, and makes the ensemble inflated inhomogeneously as well. Then, in  Section \ref{sec:l96_numeric}, we compare the performance of the baseline method and ETKF with the spectrum smoothing for various turbulence regimes and observation resolutions.

\subsubsection{Comparison with the classic ETKF \yl{in various turbulent regimes}} \label{sec:l96_numeric}

In order to demonstrate the effectiveness of the spectrum smoothing in mitigating the sampling error, \yl{we use various turbulent regimes of the Lorenz 96 models with several $F$ values (4, 8, and 16) across various ensemble sizes ($K=10, 20, 30$ and $40$) and observation resolutions (25, 33, 50, and 100\%). With a Gaussian kernel for smoothing, there are three tuning parameters: inflation level $\rho$, localization radius $c$, and smoothing kernel width $\sigma$ (the standard deviation of the Gaussian kernel). All these parameters are obtained from fine-tuning. To save computational efforts, we use a semi-joint tuning process. That is, we tune the inflation and localization first, then tune the kernel width. Once the kernel width is chosen, we tune the inflation and localization parameters again, which provides additional computational saving compared to the joint tuning of all three parameters simultaneously. 
}

The comparison of two methods for each turbulence regime is shown in Table \ref{tabl:RMSEf4} ($F=4$; weak turbulence), Table \ref{tabl:RMSEf8} ($F=8$; moderate turbulence), and Table \ref{tabl:RMSEf16} ($F=16$; \rv{strong} turbulence), respectively. For each regime, we further vary observation resolution and ensemble size to see if the relative comparison is consistent. When the turbulence is weaker, we observe more cases that the spectrum smoothing performs no better than the baseline method. However, as the turbulence becomes stronger, i.e., $F$ increases, we observe that the spectrum smoothing achieves relatively smaller (time-averaged) RMSE than the baseline method for a wider range of observation resolution and ensemble size. Although the overall RMSE increases as the system becomes more turbulent, the spectrum smoothing significantly improves the performance of baseline data assimilation. \yl{Turbulence theory is developed for a fully-developed turbulent flow with a smoother spectrum than a weak turbulent flow. Therefore, our results seem valid as smoothing the spectrum mitigates the sampling error more effectively, yielding RMSEs significantly less than those of the baseline when the turbulence is stronger.}

To compare the performance of \yl{the} two methods more easily, we shade the cases in each table when the RMSE of the baseline method is smaller than that of the spectrum smoothing. Other than those cases, the spectrum smoothing shows smaller RMSEs for all observation resolutions and ensemble sizes. As either the observation becomes sparser or the ensemble size decreases, RMSE tends to increase. However, the spectrum smoothing consistently improves the data assimilation performance compared to the baseline method when the observation resolutions \yl{are} low and the ensemble size is small. This implies that the spectrum smoothing successfully complements the information from the lack of observations or small sample (\rv{ensemble}) size using the smooth property of the turbulent flow's spectrum and, accordingly, improves the sample quality, resulting in a prominent sampling error mitigation.

\begin{table}[h]
	\centering \small
	\mbox{}\clap{
		\begin{tabular}{|c|c|c|c|c|} \hline
			\backslashbox{Observation}{K}   &10   & 20 & 30 & 40     \\ \hline
			25\%  & 0.0066 / 0.0066  &  	0.0068 / 0.0068  &	   0.0062 / 0.0062  &  	\cellcolor{Gray} 0.0061 / 0.0069 	\\ \hline	
			33\%  & 0.0073 / 0.0073  & 	    0.0061 / 0.0061   &     0.0059 / 0.0059   &  	\cellcolor{Gray} 0.0037 / 0.0058  \\ \hline	 
			50\%  & 0.0075 / 0.0075  &     0.0057 / 0.0057    &    0.0053 / 0.0053    &    	\cellcolor{Gray}  0.0041 / 0.0054 	\\ \hline	
			100\% &  0.0042 / 0.0042  &   0.0047 / 0.0047  &     0.0043 / 0.0043    &    	\cellcolor{Gray} 0.0027 / 0.0043 	\\ \hline	
		\end{tabular}
	}
	\caption{RMSE comparison between the baseline method (left column) and spectrum smoothing (right column) applied to Lorenz 96 with  $\mathbf{F=4}$ (weak turbulence) and $N=128$, $\Delta t_{obs} = 0.01 \times 15 = 0.15$, $t \in [0,200]$, observation noise=0.1854}
	\label{tabl:RMSEf4}
\end{table}

\begin{table}[!hbtp]
	\centering \small
	\mbox{}\clap{
		\begin{tabular}{|c|c|c|c|c|} \hline 
			\backslashbox{Observation}{K}   &10   & 20 & 30 & 40     \\ \hline
			25\%  & 4.5024 / 3.0801 & 4.1155 / 2.2755  & 3.8638 / 1.0186  & 3.6804 / 0.5102 \\ \hline	
			33\%  & 4.5038 / 1.0637 & 4.1080 / 0.4249 & 3.8018 / 0.2778 & 3.3590 / 0.2625 \\ \hline	 
			50\%  & 4.3720 / 0.3907  & 3.8751 / 0.2323 & 3.294 / 0.2186 &  	\cellcolor{Gray} 0.1892 / 0.1924 \\ \hline	
			100\% & 4.3437 / 0.1818  &  3.5464 / 0.1631 & 2.4276 / 0.1612 & 0.1125 / 0.1121	\\ \hline	
		\end{tabular}
	}
	\caption{RMSE comparison between the baseline method (left column) and spectrum smoothing (right column) applied to Lorenz 96 with  $\mathbf{F=8}$ (moderate turbulence) and $N=128$, $\Delta t_{obs} = 0.01 \times 15 = 0.15$, $t \in [0,200]$,  observation noise=0.3640}
	\label{tabl:RMSEf8}
\end{table}

\begin{table}[!hbtp]
	\centering \small
	\mbox{}\clap{
		\begin{tabular}{|c|c|c|c|c|} \hline
			\backslashbox{Observation}{K}   &10   & 20 & 30 & 40     \\ \hline
			25\%  & 8.0510 / 7.031  &  7.5721 / 6.9385  & 7.3357 / 6.7305  &  7.1821 / 6.4880	\\ \hline	
			33\%  & 8.0369 / 6.0142   & 7.5455 / 5.8902  &  7.2087 / 5.7011 &  6.9311 / 5.5288 	\\ \hline	 
			50\%  & 7.7317 / 2.5483 &  6.9945 / 2.3546 &  6.4049 / 1.5608   & 5.9965 / 0.6606 \\ \hline	
			100\% & 7.3057 / 0.4760  &  6.3212 / 0.3890 & 5.3280 / 0.3270 & 4.3081 / 0.3122 	\\ \hline	
		\end{tabular}
	}
	\caption{RMSE comparison \rv{between the} baseline method (left column) and spectrum smoothing (right column) applied to Lorenz 96 with  $\mathbf{F=16}$ (strong turbulence) and $N=128$, $\Delta t_{obs} = 0.01 \times 15 = 0.15$, $t \in [0,200]$, observation noise=0.6298 }
	\label{tabl:RMSEf16}
\end{table}

\subsubsection{Inhomogeneous localization and inflation effect of spectrum smoothing} \label{sec:inhomog}
\begin{figure}[h]
	\centering
	\subfloat[t=45]
	{\includegraphics[width=.4\textwidth]{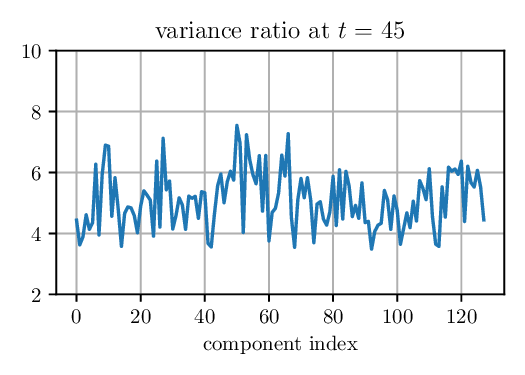}\label{fig:inhomogInf0}}	
	\subfloat[t=90]
	{\includegraphics[width=.4\textwidth]{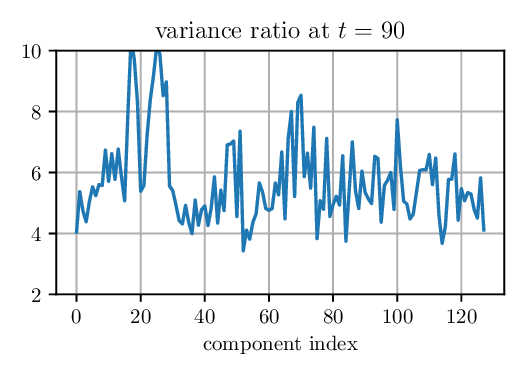}\label{fig:inhomogInf1}}\\
	\subfloat[t=135]
	{\includegraphics[width=.4\textwidth]{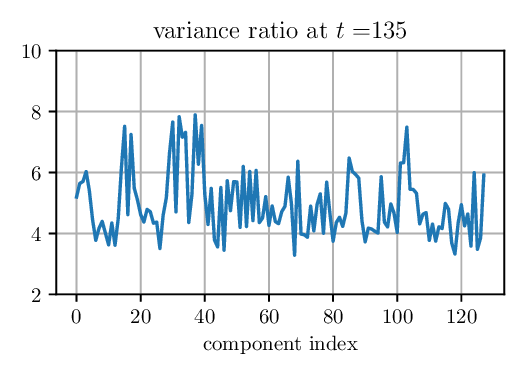}\label{fig:inhomogInf2}}		
	\subfloat[t=180]
	{\includegraphics[width=.4\textwidth]{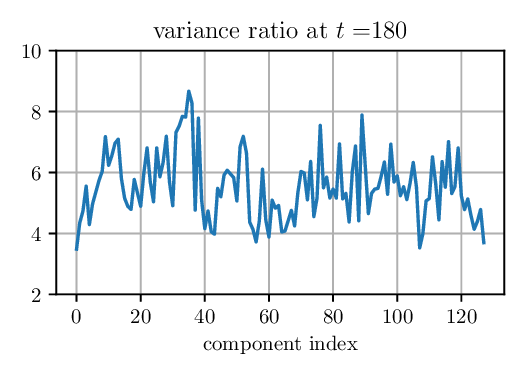}\label{fig:inhomogInf3}}\\
	\caption{\yl{Ratio of the prior variances between before and after the spectrum smoothing at various times. ETKF using Lorenz 96 model with $F=8$}}
	\label{fig:covRatio}
\end{figure}

\begin{figure}[h]
	\centering
	\subfloat[t=45]
	{\includegraphics[width=.4\textwidth]{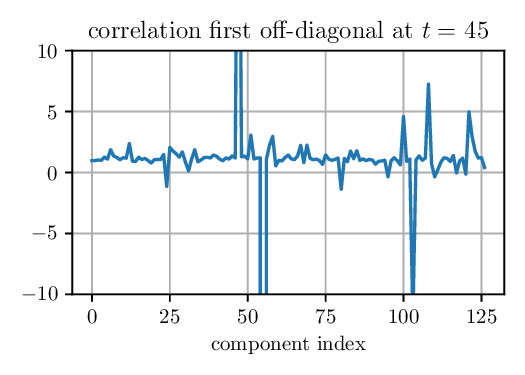}\label{fig:inhomogLoc0}}	
	\subfloat[t=90]
	{\includegraphics[width=.4\textwidth]{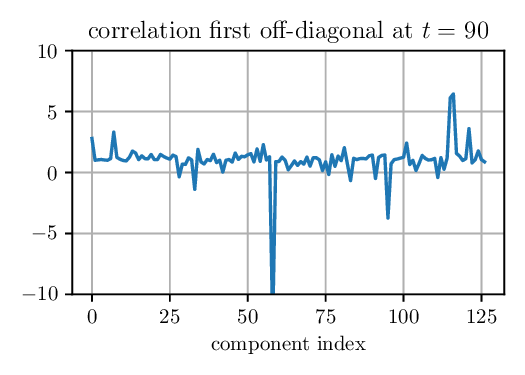}\label{fig:inhomogLoc1}}\\
	\subfloat[t=135]
	{\includegraphics[width=.4\textwidth]{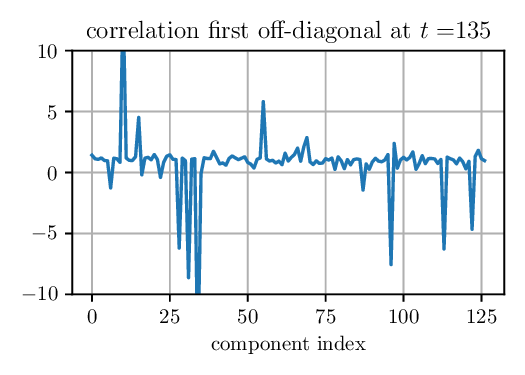}\label{fig:inhomogLoc2}}
	\subfloat[t=180]
	{\includegraphics[width=.4\textwidth]{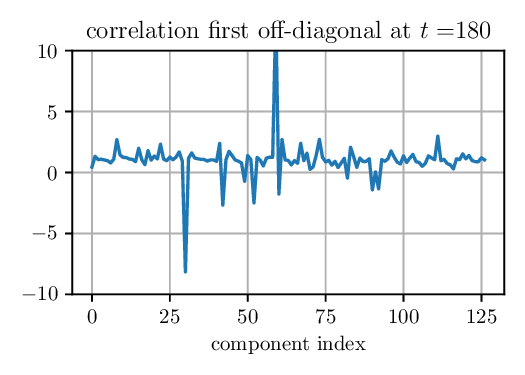}\label{fig:inhomogLoc3}}\\
	\caption{\yl{Ratio of the covariances between adjacent components (first off-diagonal components of the prior covariance matrix) between before and after the spectrum smoothing} at various times. ETKF using Lorenz 96 model with $F=8$}
	\label{fig:corrRatio}
\end{figure}

\yl{Inflation and localization are an effective way to stabilize ensemble filters and improve their performance. Due to their intrinsic nature for tuning, the effect of inflation and localization remain constant at different locations (or with different observation components). Suppose the state variable has spatially varying features. In that case, it is natural to speculate inflation and localization adaptively based on where the state variable is (in addition to the distance in the case of localization), which we call inhomogeneous inflation and localization.

To check the effect of inhomogeneous inflation and localization, we compare the diagonal and first off-diagonal components of the prior covariance i) before and ii) after the spectrum smoothing (also no multiplicative inflation and Gaspari-Cohn localization). In particular, we calculate the ratios of the variance (diagonal components) and the covariance with the adjacent component (first off-diagonal components) to check the inhomogeneous inflation and localization, which are shown in Figure \ref{fig:covRatio} and \ref{fig:corrRatio}, respectively for the Lorenz 96 with $F=8$ case. The four subplots of each Figure are at different times, $t=45, 90, 135$, and $180$. 

The ratio of variances greater than 1 in Figure \ref{fig:covRatio} implies the inflation level at each component. As seen at all four different times, the inflation level varies at different components (or locations) and is thus inhomogeneous. In particular, the inhomogeneous inflation level is not fixed in time but changes over time through the spectrum smoothing. Also, the ratio of the covariances in Figure \ref{fig:corrRatio} shows similar inhomogeneous spatially varying ratios. As this pattern does not align with the variance ratios, we can interpret that the spectrum smoothing also changes the correlation in an inhomogeneous way (that is, the correlation is modified depending on where the component is). Note that the standard Gaspari-Cohn uses a localization factor depending only on the distance between different components.
}


\section{Conclusion and Discussion} \label{sec:conclusion}

Ensemble Kalman Filter's performance is highly dependent on the ensemble size as \yl{a small ensemble size induces a sampling error and degrades the filter performance. However, due to the high computational cost of running a complex prediction model, such as a turbulent system, the typical ensemble size remains significantly small compared to the dimension of the system of interest. The current study provides another potential direction \yl{for} further investigation for mitigating sampling error through the spectrum smoothing. The spectrum smoothing is motivated by a characteristic of \yl{the} turbulent system (the Kolmogorov spectrum) and the Khinchin theorem for a stochastic process or a random field. The Kolmogorov spectrum guarantees \yl{the} smoothness of the spectrum in the dissipative regime, while the Khinchin theorem relates the spectrum with the correlation information of the (random) state. Based on these facts, the proposed method modifies the ensemble to satisfy a certain smoothness in the spectrum.

In our numerical tests using the Lorenz 96 system, the proposed spectrum smoothing shows improved filtering performance, particularly for strong turbulent regimes; as $F$ increases in the Lorenz 96 model (that is, as the model becomes more turbulent), the performance difference becomes more significant. It is worth \yl{noting} that to see performance improvement, it is essential to have a meaningful spectrum with a dissipative regime. Thus, instead of the standard 40-dimensional Lorenz 96 model, the current study considered a much larger dimension, 128, to observe performance improvements.

The current study leaves several things for further investigation. In the current study, only a Gaussian kernel has been tested for smoothing. It is natural to investigate various kernel shapes and widths for performance and sensitivity. In particular, the spectrum smoothing still requires tuning inflation, localization, and kernel width; the tuning process of three parameters can hinder the proposed method in certain application areas. In our first experiments with the Lorenz 96 model, we have an impression that the performance is less sensitive to the kernel width than the other tuning parameters (inflation and localization). Thus, a more in-depth sensitivity analysis of the spectrum smoothing (under various shapes and widths) is necessary.

Also, in the current study, we showed that (at least in the Lorenz 96 model) the spectrum smoothing is related to the effect of spatially varying inflation and localization. That is, different components of the state have different inflation and localization levels through the spectrum smoothing. As there are several analysis research works related to inflation and localization, it would be an interesting mathematical question to investigate the spectrum smoothing in the context of inflation and localization. Additionally, prior ensemble modification is a regularization using the prior knowledge of the system (which is a smooth spectrum for a turbulent system). In this line of regularization, we believe that other prior knowledge of the system in the context of regularization can also serve as a method to mitigate sampling errors, particularly related to the spatially varying regularizations \cite{inhomog_2023,inhomog2024}. 

Last, the current study tested the proposed sampling error mitigation method using only the Lorenz 96 model. However, the proposed method can be easily applied to other turbulent systems due to the non-intrusive nature of the spectrum smoothing, which we leave as future work.}

%

\subsection*{Acknowledgement}
Yoonsang Lee is supported by ONR MURI N00014-20-1-2595.

\bibliographystyle{abbrv}
\bibliography{./turbreg.bib}

\end{document}